\documentclass{article}
\usepackage{graphicx} % Required for inserting images
\usepackage{amsfonts}
\usepackage{amsmath}
\usepackage{amssymb}
\usepackage{amsthm}
\usepackage{hyperref}
\usepackage{sectsty}
\sectionfont{\centering}
\begin{document}
\title{\textbf{SOBOLEV INEQUALITIES ON K\"{A}HLER MANIFOLDS}}
\author{Sayantan Chakraborty}
\date{}
\maketitle
\begin{minipage}{0.9\textwidth}

    ABSTRACT. We prove new Sobolev type inequalities on compact K\"{a}hler manifolds with positive Ricci curvature. A proof of an already existing Sobolev inequality in the classical Bidaut-V\'{e}ron and V\'{e}ron approach is also discussed.

\end{minipage}
\section{Introduction}
We consider the class of closed Riemannian manifolds $(M,g)$ of dimension $m$ with $\text{Ric}(g) \geq (m-1)g$. For these manifolds, we have the following Sobolev inequality
\begin{align*}
    \Bigg(\frac{1}{\text{Vol}(M)}\int_{M}|\phi|^{q+1}\Bigg)^{\frac{2}{q+1}} \leq \frac{q-1}{m}\Bigg(\frac{1}{\text{Vol}(M)}\int_{M}|\nabla \phi|^2\Bigg)+\frac{1}{\text{Vol}(M)}\int_{M}\phi^2
\end{align*}
for $1<q \leq \frac{m+2}{m-2}$. This was proved by Bidaut-V\'{e}ron and V\'{e}ron in [1]. \newline\newline
The constant $\frac{q-1}{m}$ cannot be improved in this class of Riemannian manifolds since it is optimal for spheres (see for instance [2]).\newline\newline
K\"{a}hler manifolds are a special class of Riemannian manifolds with three mutually compatible structures: a complex structure, a Riemannian structure and a symplectic structure. Due to this special feature, improvements of some results for Riemannian manifolds are possible in the K\"{a}hler setting. We refer interested readers to [3, 7, 8, 4, 5, 9, 10] and references therein for interesting comparison theorems involving eigenvalues, volume and other geometric quantities. \newline\newline
Baudoin and Munteanu in their paper [6] improved the constant $\frac{q-1}{m}$ in the K\"{a}hler setting for the range $1<q<\frac{m+2}{m-2}$ where $m=\text{dim}_{\mathbb{R}}M$. They showed the following:\newline\newline
\textbf{Theorem 1.} (Baudoin-Munteanu [6]). \textit{Let} $(M,\omega)$ \textit{be a compact K\"{a}hler manifold of complex dimension} $n$ \textit{with} $\text{Ric}(\omega) \geq \omega$. \textit{For any} $1<q<\frac{n+1}{n-1}$ \textit{we have the Sobolev inequality}
\begin{align*}
    \Bigg(\frac{1}{\text{Vol}(M)}\int_{M} |\phi|^{q+1}\Bigg)^{\frac{2}{q+1}}\leq \frac{C_S}{\text{Vol}(M)}\int_{M}|\nabla \phi|^2+\frac{1}{\text{Vol}(M)}\int_{M}\phi^2
\end{align*}
\textit{for any} $\phi \in C^{\infty}(M)$ \textit{where}
\begin{align*}
    C_S=(q-1)\Bigg(\frac{2n+q+2-2\sqrt{(n+1)(n+1-(n-1)q)}}{2qn}\Bigg).
\end{align*}
\newline
Let us set the notation used throughout this note. Suppose $(M,\omega)$ is a compact K\"{a}hler manifold of complex dimension $n$ with complex structure $J$. If $\{e_k\}_{k=1}^{2n}$ is an orthonormal frame such that $e_{2i}=J(e_{2i-1})$ for $i \in \{1,2,\cdots,n\}$ then we work in the following unitary frame
\begin{align*}
    X_{i}=\frac{1}{\sqrt{2}}(e_{2i-1}-\sqrt{-1}e_{2i}), \,X_{\bar{i}}=\frac{1}{\sqrt{2}}(e_{2i-1}+\sqrt{-1}e_{2i})
\end{align*}
for $i \in \{1,2,\cdots,n\}$. These satisfy $\{X_{k},X_{\bar{l}}\}=\delta_{kl}$ and $\{X_{k},X_{l}\}=0$. The Ricci curvature condition is then equivalent to $R_{i\bar{j}} \geq \delta_{i\bar{j}}$. The complex Laplace operator is given by $\Box u=\sum_{i=1}^{n}u_{i\bar{i}}$. The real Laplace operator $\Delta$ relates to this as $\Delta=2\Box$. We also have $|\partial u|^2=\sum_{i=1}^{n}u_{i}u_{\bar{i}}=\frac{1}{2}|\nabla u|^2$.\newline\newline
\textbf{Remark:} Let $f$ be an eigenfunction for the first nonzero Dirichlet eigenvalue $\lambda_1$ of the complex Laplace operator $\Box$ on $M$, i.e $-\Box f=\lambda_1 f$. Let us set $\phi=1+tf$ in theorem 1. We aim to compare the coefficient of $t^2$ on both sides of the inequality. Let
\begin{align*}
    \psi(t)=\Bigg(\frac{1}{\text{Vol}(M)}\int_{M}|1+tf|^{q+1}\Bigg)^{\frac{2}{q+1}}.
\end{align*}
Then near $t=0$ we have
\begin{align*}
    \psi'(t)=\frac{2}{\text{Vol}(M)}\Bigg(\int_{M}(1+tf)^{q}f\Bigg)\Bigg(\frac{1}{\text{Vol}(M)}\int_{M}(1+tf)^{q+1}\Bigg)^{\frac{2}{q+1}-1}
\end{align*}
and
\begin{align*}
    &\psi''(t)=\frac{2}{\text{Vol}(M)}\Bigg(\frac{1}{\text{Vol}(M)}\int_{M}(1+tf)^{q+1}\Bigg)^{\frac{2}{q+1}-2}\\
    &\,\,\,\,\,\,\,\,\,\,\,\,\,\,\,\Bigg((1-q)\frac{1}{\text{Vol}(M)}\Big(\int_{M}(1+tf)^{q}f\Big)^2+\Big(q\int_{M}(1+tf)^{q-1}f^2\Big)\Bigg(\frac{1}{\text{Vol}(M)}\int_{M}(1+tf)^{q+1}\Bigg)\Bigg).\\
\end{align*}
In particular the coefficient of $t^2$ on the left hand side of the inequality is
\begin{align*}
    \frac{1}{2}\psi''(0)=\frac{q}{\text{Vol}(M)}\int_{M}f^2
\end{align*}
where we use the fact that $\int_{M}f=-\frac{1}{\lambda_1}\int_{M}\Box f=0$.
\newline\newline
The right hand side of the inequality with $\phi=1+tf$ is
\begin{align*}
    & \,\,\,\,\,\,\,\,\frac{t^2C_S}{\text{Vol}(M)}\int_{M}|\nabla f|^2+\frac{1}{\text{Vol}(M)}\int_{M}(1+tf)^2\\
    & =\frac{2t^2C_S}{\text{Vol}(M)}\int_{M}|\partial f|^2+\frac{1}{\text{Vol}(M)}\int_{M}(1+tf)^2.
\end{align*}
We have $\int_{M}|\partial f|^2=-\int_{M}f\Box f=\lambda_1\int_{M}f^2$. Therefore the coeffient of $t^2$ on the right hand side is
\begin{align*}
    \frac{2C_S\lambda_1+1}{\text{Vol}(M)}\int_{M}f^2.
\end{align*}
Comparing these coeffients we get
\begin{align*}
    2C_S\lambda_1+1 \geq q \Rightarrow \lambda_1 \geq \frac{q-1}{2C_S}.
\end{align*}
On the class of K\"{a}hler manifolds with $\text{Ric}(\omega) \geq \omega$, we have an improved Lichnerowicz estimate $\lambda_1 \geq 1$ (see for example Futaki [11, Theorem 2.4.5]). For this reason we conjecture the following:\newline\newline
\textbf{Conjecture:} $C_S$ can be improved to $\frac{q-1}{2}$ in theorem 1. \newline\newline
There are a lot of K\"{a}hler manifolds with $\text{Ric}(\omega) \geq \omega$ and $\lambda_1=1$. Some examples of such manifolds are \begin{align*}
    (M,\omega)=(\mathbb{CP}^{k},\omega_{\mathbb{CP}^{k}}) \times (\mathbb{CP}^{l}, (1-a)\omega_{\mathbb{CP}^{l}}).
\end{align*}
where $k$, $l$ are positive integers with $k+l=n$ and $a \in [0,1)$. In fact, any K\"{a}hler-Einstein Fano manifold admitting a nontrivial holomorphic vector field has $\lambda_1=1$ [12]. This means that $\frac{q-1}{2}$ is the best possible candidate for $C_S$.\newline\newline
 Theorem 1 can be proven using the standard approach of working with an associated PDE. Namely, we consider the functional
\begin{align*}
    \mathcal{F}(u)=\frac{\int_{M}|\partial u|^2+\lambda\int_{M}u^2}{\Big(\int_{M}|u|^{q+1}\Big)^{\frac{2}{q+1}}},\,u\in W^{1,2}(M)\backslash\{0\}.
\end{align*}
The first variation in the direction of $\dot{u}$ is
\begin{align*}
& 2\Bigg[\frac{\int_{M}\langle\partial u,\partial \dot{u}\rangle+\lambda\int_{M}u\dot{u}}{\Big(\int_{M}|u|^{q+1}\Big)^{\frac{2}{q+1}}}-\frac{\int_{M}|u|^q\dot{u}}{\Big(\int_{M}|u|^{q+1}\Big)^{\frac{2}{q+1}+1}}\Big(\int_{M}|\partial u|^2+\lambda \int_{M} u^2\Big)\Bigg]\\
&=\frac{2}{\Big(\int_{M}|u|^{q+1}\Big)^{\frac{2}{q+1}}}\Bigg[\int_{M}(-\dot{u}\Box{u}+\lambda u\dot{u})-\frac{\int_{M}|\partial u|^2+\lambda\int_{M}u^2}{\int_{M}|u|^{q+1}}\int_{M}|u|^q\dot{u}\Bigg].
\end{align*}
Thus a positive $u$ is a critical point of $\mathcal{F}$ iff
\begin{align*}
    -\Box u+\lambda u=cu^q
\end{align*}
on $M$, with $c=\frac{\int_{M}|\partial u|^2+\lambda\int_{M}u^2}{\int_{M} |u|^{q+1}}$. The unknown function can be rescaled so that $c=1$.\newline\newline
Therefore theorem 1 follows from the next theorem.\newline\newline
\textbf{Theorem A.} (Baudoin-Munteanu [6]). \textit{Let} $(M,\omega)$ \textit{be a compact K\"{a}hler manifold of complex dimension} $n$ \textit{and} $\text{Ric}(\omega) \geq \omega$. \textit{If} $0<\lambda<\frac{1}{2C_S}$\textit{, then the only positive solutions of}
\begin{align*}
    -\Box u+\lambda u=u^{q}
\end{align*}
\textit{are positive constant functions.}\newline\newline
We plan to have an alternate proof of theorem A in the next section. Our style of proving it is in the flavour of the Bidaut-V\'{e}ron and V\'{e}ron approach.\newline\newline
\textbf{Remark:} While discussing the proof of theorem A in the next section, we explain at the end of lemma 2.3 why we fail to get $C_S=\frac{q-1}{2}$ using the coefficients obtained there. The coefficients also point us to the fact that just using an integral identity using the complex Hessian doesn't help get a better Sobolev constant (as opposed to the proof of $\lambda_1\geq 1$ for K\"{a}hler manifolds with $\text{Ric}(\omega) \geq \omega$ in [4] for example, which uses an integral involving only the complex Hessian). Therefore we need some support from an integral identity (lemma 4 in the next section) using the other part of the Hessian.\newline\newline
In this note, we prove a similar PDE statement where the constant $\lambda$ involves the first Dirichlet eigenvalue $\lambda_1$ of the complex Laplace operator $\Box$ on a compact K\"{a}hler manifold with positive Ricci lower bound. The following is the main result of the paper:\newline\newline
\textbf{Theorem B.} \textit{Let} $(M,\omega)$ \textit{be a compact K\"{a}hler manifold of complex dimension} $n$ \textit{and} $\text{Ric}(\omega) \geq \omega$. \textit{If} $0<\lambda<\frac{1}{2C_{S,\lambda_1,k}}$\textit{, then the only positive solutions of}
\begin{align*}
    -\Box u+\lambda u=u^q
\end{align*}
\textit{are positive constant functions. Here the constant} $C_{S,\lambda_1,k}$ \textit{is given by}
\begin{align*}
    \frac{1}{C_{S,\lambda_1,k}}=\frac{2}{q-1}\Bigg(\Bigg(1-\frac{(n+(n-1)k)(kn+n-1)q}{(4n^2+4n+q)k}\Bigg)\lambda_1+\frac{qn(kn+n-1)}{(4n^2+4n+q)k}\Bigg)
\end{align*}
\textit{where} $k$ \textit{is any number lying in the interval}
\begin{align*}
    \Bigg[\frac{2(n+1)}{q(n-1)}-\frac{2(n+1)}{q(n-1)}\sqrt{1-\frac{(n-1)q}{n+1}}-1,\frac{2(n+1)}{q(n-1)}+\frac{2(n+1)}{q(n-1)}\sqrt{1-\frac{(n-1)q}{n+1}}-1\Bigg].
    \end{align*}
    We note that the coeffient of $\lambda_1$ in the expression for $\frac{1}{C_{S,\lambda_1,k}}$ is nonnegative if $k$ lies in the above interval.
\newline\newline
As a corollary we obtain the following Sobolev inequality:\newline\newline
\textbf{Theorem 2.} \textit{Let} $(M,\omega)$ \textit{be a compact K\"{a}hler manifold of complex dimension} $n$ \textit{and} $\text{Ric}(\omega)\geq \omega$. \textit{For any} $1 \leq q<\frac{n+1}{n-1}$ \textit{we have the Sobolev inequality}
\begin{align*}
    \Bigg(\frac{1}{\text{Vol}(M)}\int_{M}|\phi|^{q+1}\Bigg)^{\frac{2}{q+1}} \leq \frac{1}{\text{Vol}(M)}\int_{M}\phi^2+\frac{C_{S,\lambda_1,k}}{\text{Vol}(M)}\int_{M}|\nabla \phi|^2
\end{align*}
\textit{for any} $\phi \in C^{\infty}(M)$. \textit{Here} $C_{S,\lambda_1,k}$ \textit{is as in} theorem B \textit{and} $k$ \textit{is any number lying in the interval}
\begin{align*}
    \Bigg[\frac{2(n+1)}{q(n-1)}-\frac{2(n+1)}{q(n-1)}\sqrt{1-\frac{(n-1)q}{n+1}}-1,\frac{2(n+1)}{q(n-1)}+\frac{2(n+1)}{q(n-1)}\sqrt{1-\frac{(n-1)q}{n+1}}-1\Bigg].
    \end{align*}\newline\newline
\textbf{Remark:} When $k=\frac{2(n+1)}{q(n-1)}-\frac{2(n+1)}{q(n-1)}\sqrt{1-\frac{(n-1)q}{n+1}}-1$, the coefficient of $\lambda_1$ vanishes and we get $C_{S,\lambda_1,k} =C_S$ where $C_S$ is the known Sobolev constant for K\"{a}hler manifolds with $\text{Ric}(\omega) \geq \omega$. So this inequality can be viewed as a generalized Sobolev inequality involving $\lambda_1$.\newline\newline
\textbf{Acknowledgement:} It is my pleasure to thank my advisor Professor Xiaodong Wang for introducing me to this problem and for many helpful discussions.
\section{Proof of theorem A}
Let $u \in C^{\infty}(M)$ be a positive solution of
\begin{align*}
    -\Box u+\lambda u=u^q.
\end{align*}
We set $u=v^{-\beta}$, $v>0$, $\beta \in \mathbb{R}-\{0\}$. Then $v$ satisfies the equation
\begin{align*}
    \Box v=\frac{1}{\beta}v^{\beta+1-\beta q}-\frac{\lambda}{\beta}v+(\beta+1)\frac{|\partial v|^2}{v}.\tag{*}
\end{align*}
For calculations we omit the summation notation as this will always be clear. Theorem A will be proven in steps. Here is our first lemma.\newline\newline\newline
\textbf{Lemma 2.1}: Let $(M,\omega)$ be a compact K\"{a}hler manifold satisfying $\text{Ric}(\omega) \geq \omega$. For $\gamma \geq 0$ and a positive $v \in C^{\infty}(M)$, we have
\begin{align*}
    &\int_{M}v^{\gamma}|v_{\bar{i}\bar{j}}|^2 \leq \gamma(\gamma-1)\int_{M}v^{\gamma-2}|\partial v|^4+\gamma\int_{M}v^{\gamma-1}v_{\bar{i}}v_{j}v_{i\bar{j}}+2\gamma\int_{M}v^{\gamma-1}|\partial v|^2\Box v\\
    &\,\,\,\,\,\,\,\,\,\,\,\,\,\,\,\,\,\,\,\,\,\,\,\,\,\,\,\,\,\,\,\,\,\,\,\,\,+\int_{M}v^{\gamma}(\Box v)^2-\int_{M}v^{\gamma}|\partial v|^2.\,\,\,\,\,\,\,\,\,\,\,\,\,\,\,\,\,\,\,\,\,\,\,\,\\
\end{align*}
\begin{proof}
 We compute
\begin{align*}
&\,\,\,\,\,\,\,\,\,\,\,\,\, \int_{M} v^{\gamma}|v_{\bar{i}\bar{j}}|^2 =\int_{M}v^{\gamma}v_{\bar{i}\bar{j}}v_{ij}=-\int_{M}(v^{\gamma}v_{\bar{i}\bar{j}})_{j}v_{i}\\
    &=-\gamma\int_{M}v^{\gamma-1}v_{i}v_{j}v_{\bar{i}\bar{j}}-\int_{M}v^{\gamma}v_{\bar{i}\bar{j},j}v_{i}\\
    &=-\gamma\int_{M}v^{\gamma-1}v_{i}v_{j}v_{\bar{i}\bar{j}}+\int_{M}v^{\gamma}(-v_{j\bar{j},\bar{i}}+R_{j\bar{i}l\bar{j}}v_{\bar{l}})v_{i}\\
    &=-\gamma\int_{M}v^{\gamma-1}v_{\bar{i}\bar{j}}v_{i}v_{j}-\int_{M}v^{\gamma}(\Box v)_{\bar{i}}v_{i}-\int_{M}v^{\gamma}R_{l\bar{i}}v_{i}v_{\bar{l}}\\
    & \leq -\gamma\int_{M}v^{\gamma-1}v_{\bar{i}\bar{j}}v_{i}v_{j}+\int_{M} v^{\gamma}(\Box v)^2+\gamma\int_{M}v^{\gamma-1}|\partial v|^2\Box v-\int_{M}v^{\gamma}|\partial v|^2\\
    &=\gamma\int_{M} (v^{\gamma-1}v_{i}v_{j})_{\bar{j}}v_{\bar{i}}+\int_{M} v^{\gamma}(\Box v)^2+\gamma\int_{M}v^{\gamma-1}|\partial v|^2\Box v-\int_{M}v^{\gamma}|\partial v|^2\\
    &=\gamma(\gamma-1)\int_{M}v^{\gamma-2}v_{\bar{j}}v_{i}v_{j}v_{\bar{i}}+\gamma\int_{M}v^{\gamma-1}v_{i\bar{j}}v_{j}v_{\bar{i}}+\gamma\int_{M}v^{\gamma-1}v_{i}v_{j\bar{j}}v_{\bar{i}}\\
    & \,\,\,\,\,+\int_{M} v^{\gamma}(\Box v)^2+\gamma\int_{M}v^{\gamma-1}|\partial v|^2\Box v-\int_{M}v^{\gamma}|\partial v|^2\\
    &=\gamma(\gamma-1)\int_{M}v^{\gamma-2}|\partial v|^4+\gamma\int_{M}v^{\gamma-1}v_{\bar{i}}v_{j}v_{i\bar{j}}+2\gamma\int_{M}v^{\gamma-1}|\partial v|^2\Box v\\
    &\,\,\,\,\,+\int_{M}v^{\gamma}(\Box v)^2-\int_{M}v^{\gamma}|\partial v|^2.
\end{align*}
Here in the process we have integrated by parts repeatedly. We have used the Ricci identities while going from the 4th to the 5th expression and the Ricci curvature condition while going from the 6th to the 7th expression.
\end{proof}
We plan to write an inequality that uses lemma 2.1 and introduces an additional parameter. Before doing that we need some calculations.\newline\newline
\textbf{Lemma 2.2}: Suppose $v$ is the function satisfying (*). Then for $\gamma>0$ we have
\begin{align*}
    &\int_{M}v^{\gamma-1}|\partial v|^2\Box v=\frac{1}{\beta}\int_{M}v^{\beta+\gamma-\beta q}|\partial v|^2-\frac{\lambda}{\beta}\int_{M}v^{\gamma}|\partial v|^2+(\beta+1)\int_{M}v^{\gamma-2}|\partial v|^4. \tag 1\\
    &\int_{M}v^{\gamma}(\Box v)^2=\frac{\beta q-\gamma}{\beta}\int_{M}v^{\beta+\gamma-\beta q}|\partial v|^2+\frac{\lambda(\gamma-\beta)}{\beta}\int_{M}v^{\gamma}|\partial v|^2\\
    &\,\,\,\,\,\,\,\,\,\,\,\,\,\,\,\,\,\,\,\,\,\,\,\,\,\,\,\,\,\,\,\,\,\,\,\,\,\,\,+(\beta+1)^2\int_{M}v^{\gamma-2}|\partial v|^4.\tag 2\\
    &\int_{M}v^{\gamma-1}v_{i\bar{j}}v_{\bar{i}}v_{j}=\frac{1}{\gamma}\int_{M}(\Box v)^2v^{\gamma}+\int_{M}v^{\gamma-1}|\partial v|^2\Box v-\frac{1}{\gamma}\int_{M}|v_{i\bar{j}}|^2v^{\gamma}.\tag 3\\
\end{align*}
\begin{proof}
    Using the equation satisfied by $\Box v$ we have
    \begin{align*}
    &\,\,\,\,\,\,\, \int_{M}v^{\gamma-1}|\partial v|^2\Box v\\
    &=\int_{M}v^{\gamma-1}|\partial v|^2\Big(\frac{1}{\beta}v^{\beta+1-\beta q}-\frac{\lambda}{\beta}v+(\beta+1)\frac{|\partial v|^2}{v}\Big)\\
    &=\frac{1}{\beta}\int_{M}v^{\beta+\gamma-\beta q}|\partial v|^2-\frac{\lambda}{\beta}\int_{M}v^{\gamma}|\partial v|^2+(\beta+1)\int_{M}v^{\gamma-2}|\partial v|^4.
\end{align*}
This proves (1). For the second equality we compute
    \begin{align*}
    &\,\,\,\,\,\,\,\,\int_{M} v^{\gamma}(\Box v)^2\\
    &=\int_{M}v^{\gamma}\Box v\Big(\frac{1}{\beta}v^{\beta+1-\beta q}-\frac{\lambda}{\beta}v+(\beta+1)\frac{|\partial v|^2}{v}\Big)\\
    &=\frac{1}{\beta}\int_{M}v^{\beta+\gamma+1-\beta q}\Box v-\frac{\lambda}{\beta}\int_{M}v^{\gamma+1}\Box v+(\beta+1)\int_{M}v^{\gamma-1}\Box v|\partial v|^2\\
    &=\frac{\beta q-\beta-\gamma-1}{\beta}\int_{M}v^{\beta+\gamma-\beta q}|\partial v|^2+\frac{\lambda(\gamma+1)}{\beta}\int_{M}v^{\gamma}|\partial v|^2\\
    &\,\,\,\,\,\,+(\beta+1)\int_{M}v^{\gamma-1}\Box v|\partial v|^2.\\
\end{align*}
Using the first equality we get
\begin{align*}
    & \,\,\,\,\,\int_{M}v^{\gamma}(\Box v)^2\\
    &=\frac{\beta q-\beta-\gamma-1}{\beta}\int_{M}v^{\beta+\gamma-\beta q}|\partial v|^2+\frac{\lambda(\gamma+1)}{\beta}\int_{M}v^{\gamma}|\partial v|^2\\
    &\,\,\,\,\,\,\,\,\,+\frac{(\beta+1)}{\beta}\int_{M}v^{\beta+\gamma-\beta q}|\partial v|^2-\frac{(\beta+1)\lambda}{\beta}\int_{M}v^{\gamma}|\partial v|^2+(\beta+1)^2\int_{M}v^{\gamma-2}|\partial v|^4\\
    &=\frac{\beta q-\gamma}{\beta}\int_{M}v^{\beta+\gamma-\beta q}|\partial v|^2+\frac{\lambda(\gamma-\beta)}{\beta}\int_{M}v^{\gamma}|\partial v|^2+(\beta+1)^2\int_{M}v^{\gamma-2}|\partial v|^4.
\end{align*}
This proves (2).\newline\newline
For (3) we compute
\begin{align*}
    & \,\,\,\,\,\,\,\int_{M}|v_{i\bar{j}}|^2v^{\gamma}\\
    &=\int_{M}v_{i\bar{j}}v_{j\bar{i}}v^{\gamma}\\
    &=-\int_{M}(v_{i\bar{j}}v^{\gamma})_{j}v_{\bar{i}}\\
    &=-\int_{M}v_{i\bar{j},j}v^{\gamma}v_{\bar{i}}-\gamma\int_{M}v_{i\bar{j}}v^{\gamma-1}v_{j}v_{\bar{i}}\\
    &=-\int_{M}v_{j\bar{j},i}v^{\gamma}v_{\bar{i}}-\gamma\int_{M}v_{i\bar{j}}v^{\gamma-1}v_{\bar{i}}v_{j}\\
    &=-\int_{M}(\Box v)_{i}v_{\bar{i}}v^{\gamma}-\gamma\int_{M}v_{i\bar{j}}v^{\gamma-1}v_{\bar{i}}v_{j}\\
    &=\int_{M} (\Box v)^2v^{\gamma}+\gamma\int_{M}(\Box v)v_{\bar{i}}v^{\gamma-1}v_{i}-\gamma\int_{M}v_{i\bar{j}}v^{\gamma-1}v_{\bar{i}}v_{j}\\
    &=\int_{M} (\Box v)^2v^{\gamma}+\gamma\int_{M}v^{\gamma-1}|\partial v|^2(\Box v)-\gamma\int_{M}v^{\gamma-1}v_{i\bar{j}}v_{\bar{i}}v_{j}.\\
\end{align*}
Therefore
\begin{align*}
    \int_{M}v^{\gamma-1}v_{i\bar{j}}v_{\bar{i}}v_{j}=\frac{1}{\gamma}\int_{M}(\Box v)^2v^{\gamma}+\int_{M}v^{\gamma-1}|\partial v|^2\Box v-\frac{1}{\gamma}\int_{M}|v_{i\bar{j}}|^2v^{\gamma}.
\end{align*}
\end{proof}
Here is the lemma involving the additional parameter.\newline\newline
\textbf{Lemma 2.3} For the function $v$ satisfying (*) and $\gamma>0, a \in \mathbb{R}$ we have
\begin{align*}
&\int_{M}\Big|v_{\bar{i}\bar{j}}+a\frac{v_{\bar{i}}v_{\bar{j}}}{v}\Big|^2v^{\gamma} \leq A_1\int_{M}v^{\gamma-2}|\partial v|^4 + B_1\int_{M}v^{\beta+\gamma-\beta q}|\partial v|^2\\
&\,\,\,\,\,\,\,\,\,\,\,\,\,\,\,\,\,\,\,\,\,\,\,\,\,\,\,\,\,\,\,\,\,\,\,\,\,\,\,\,\,\,\,\,\,\,\,\,\,\,\,\,\,\,\,\,\,\,\,\,\,+C_1\int_{M}v^{\gamma}|\partial v|^2+D_1\int_{M}|v_{i\bar{j}}|^2v^{\gamma}
\end{align*}
where
\begin{align*}
    & A_1=\gamma(\gamma-1)-2a(\gamma-1)+a^2+(3\gamma-4a)(\beta+1)+\Big(2-\frac{2a}{\gamma}\Big)(\beta+1)^2\\
    & B_1=\frac{3\gamma-4a}{\beta}+\Big(2-\frac{2a}{\gamma}\Big)\Big(\frac{\beta q-\gamma}{\beta}\Big)\\
    & C_1=\frac{(4a-3\gamma)\lambda}{\beta}+\Big(2-\frac{2a}{\gamma}\Big)\frac{\lambda(\gamma-\beta)}{\beta}-1\\
    & D_1=\frac{2a}{\gamma}-1.\\
\end{align*}
\begin{proof}
    \begin{align*}
    &\int_{M}\Big|v_{\bar{i}\bar{j}}+a\frac{v_{\bar{i}}v_{\bar{j}}}{v}\Big|^2v^{\gamma}\\
    &=\int_{M}v^{\gamma}|v_{\bar{i}\bar{j}}|^2+a\int_{M}v_{\bar{i}\bar{j}}v_{i}v_{j}v^{\gamma-1}+a\int_{M}v_{ij}v_{\bar{i}}v_{\bar{j}}v^{\gamma-1}+a^2\int_{M}v^{\gamma-2}|\partial v|^4. \tag 4
\end{align*}
We have
\begin{align*}
    & \int_{M}v^{\gamma-1}v_{i}v_{j}v_{\bar{i}\bar{j}}\\
    &=-\int_{M}(v^{\gamma-1}v_{i}v_{j})_{\bar{j}}v_{\bar{i}}\\
    &=-\int_{M}v^{\gamma-1}v_{i\bar{j}}v_{j}v_{\bar{i}}-(\gamma-1)\int_{M}v^{\gamma-2}v_{\bar{j}}v_{i}v_{j}v_{\bar{i}}-\int_{M}v^{\gamma-1}v_{i}\Box v v_{\bar{i}}\\
    &=-\int_{M}v^{\gamma-1}v_{i\bar{j}}v_{\bar{i}}v_{j}-(\gamma-1)\int_{M}v^{\gamma-2}|\partial v|^4-\int_{M}v^{\gamma-1}|\partial v|^2\Box v.
\end{align*}
We note that all the integrals appearing in the last line are real numbers. Therefore the second and third terms in the right hand side of (4) are real and equal. Hence (4) becomes
\begin{align*}
    &\int_{M}v^{\gamma}|v_{\bar{i}\bar{j}}|^2-2a\int_{M}v^{\gamma-1}v_{i\bar{j}}v_{\bar{i}}v_{j}-2a(\gamma-1)\int_{M}v^{\gamma-2}|\partial v|^4\\
    &-2a\int_{M}v^{\gamma-1}|\partial v|^2\Box v +a^2\int_{M}v^{\gamma-2}|\partial v|^4.
\end{align*}
By lemma 2.1 and part (3) of lemma 2.2 this expression is less than or equal to
\begin{align*}
& \gamma(\gamma-1)\int_{M}v^{\gamma-2}|\partial v|^4+\gamma\int_{M}v^{\gamma-1}v_{\bar{i}}v_{j}v_{i\bar{j}}+2\gamma\int_{M}v^{\gamma-1}|\partial v|^2\Box v+\int_{M}v^{\gamma}(\Box v)^2\\
&\,\,\,\,\,-\int_{M}v^{\gamma}|\partial v|^2-2a\int_{M}v^{\gamma-1}v_{i\bar{j}}v_{\bar{i}}v_{j}-2a(\gamma-1)\int_{M}v^{\gamma-2}|\partial v|^4\\
&\,\,\,\,\,-2a\int_{M}v^{\gamma-1}|\partial v|^2\Box v +a^2\int_{M}v^{\gamma-2}|\partial v|^4\\
&=\gamma(\gamma-1)\int_{M}v^{\gamma-2}|\partial v|^4+(\gamma-2a)\Bigg(\frac{1}{\gamma}\int_{M}(\Box v)^2v^{\gamma}+\int_{M}v^{\gamma-1}|\partial v|^2\Box v-\frac{1}{\gamma}\int_{M}|v_{i\bar{j}}|^2v^{\gamma}\Bigg)\\
&\,\,\,\,\,+2(\gamma-a)\int_{M}v^{\gamma-1}|\partial v|^2\Box v +(a^2-2a(\gamma-1))\int_{M}v^{\gamma-2}|\partial v|^4 +\int_{M}v^{\gamma}(\Box v)^2-\int_{M}v^{\gamma}|\partial v|^2\\
    &=(\gamma(\gamma-1)-2a(\gamma-1)+a^2)\int_{M}v^{\gamma-2}|\partial v|^4+(3\gamma-4a)\int_{M}v^{\gamma-1}|\partial v|^2\Box v\\
    &\,\,\,\,\,\,\,+\Big(2-\frac{2a}{\gamma}\Big)\int_{M}(\Box v)^2v^{\gamma}+\Big(\frac{2a}{\gamma}-1\Big)\int_{M}|v_{i\bar{j}}|^2v^{\gamma}-\int_{M}v^{\gamma}|\partial v|^2.
\end{align*}
Using equations (1) and (2) of lemma 2.2, this expression becomes
\begin{align*}
    &\,\,\,\,\,\,\,(\gamma(\gamma-1)-2a(\gamma-1)+a^2)\int_{M}v^{\gamma-2}|\partial v|^4\\
    &+(3\gamma-4a)\Big(\frac{1}{\beta}\int_{M}v^{\beta+\gamma-\beta q}|\partial v|^2-\frac{\lambda}{\beta}\int_{M}v^{\gamma}|\partial v|^2+(\beta+1)\int_{M}v^{\gamma-2}|\partial v|^4\Big)\\
    &+\Big(2-\frac{2a}{\gamma}\Big)\Bigg(\frac{\beta q-\gamma}{\beta}\int_{M}v^{\beta+\gamma-\beta q}|\partial v|^2+\frac{\lambda(\gamma-\beta)}{\beta}\int_{M}v^{\gamma}|\partial v|^2+(\beta+1)^2\int_{M}v^{\gamma-2}|\partial v|^4\Bigg)\\
    &+\Big(\frac{2a}{\gamma}-1\Big)\int_{M}|v_{i\bar{j}}|^2v^{\gamma}-\int_{M}v^{\gamma}|\partial v|^2.
\end{align*}
Collecting like terms together we get
\begin{align*}
    & \Bigg(\gamma(\gamma-1)-2a(\gamma-1)+a^2+(3\gamma-4a)(\beta+1)+\Big(2-\frac{2a}{\gamma}\Big)(\beta+1)^2\Bigg)\int_{M}v^{\gamma-2}|\partial v|^4\\
    &+\Bigg(\frac{3\gamma-4a}{\beta}+\Big(2-\frac{2a}{\gamma}\Big)\Big(\frac{\beta q-\gamma}{\beta}\Big)\Bigg)\int_{M}v^{\beta+\gamma-\beta q}|\partial v|^2\\
    &+\Bigg(\frac{(4a-3\gamma)\lambda}{\beta}+\Big(2-\frac{2a}{\gamma}\Big)\frac{\lambda(\gamma-\beta)}{\beta}-1\Bigg)\int_{M}v^{\gamma}|\partial v|^2+\Big(\frac{2a}{\gamma}-1\Big)\int_{M}|v_{i\bar{j}}|^2v^{\gamma}.
\end{align*}
\end{proof}
\textit{Revisiting the remark at the end of theorem A from the previous section}:\newline
We aim to choose the parameters in lemma 2.3 in a way so as to make the coefficients $A_1,B_1,C_1,D_1$ nonpositive with atleast one of them being strictly negative. This will force $v$ to be a constant. The conditions $B_1 \leq 0$ and $C_1 \leq 0$ combined give us
\begin{align*}
    &\Big(2-\frac{2a}{\gamma}\Big)\Big(q-\frac{\gamma}{\beta}\Big)\lambda+\Big(2-\frac{2a}{\gamma}\Big)\lambda\Big(\frac{\gamma}{\beta}-1\Big)-1\\
    & \leq \Big(\frac{4a-3\gamma}{\beta}\Big)\lambda+\Big(2-\frac{2a}{\gamma}\Big)\frac{\lambda(\gamma-\beta)}{\beta}-1 \leq 0
\end{align*}
i.e,
\begin{align*}
    \Big(2-\frac{2a}{\gamma}\Big)(q-1)\lambda-1 \leq 0.
\end{align*}
The condition $D_1 \leq 0$ is equivalent to
\begin{align*}
    \Big(2-\frac{2a}{\gamma}\Big) \geq 1>0
\end{align*}
The two inequalities combined give us the condition
\begin{align*}
    \lambda \leq \frac{1}{(q-1)\Big(2-\frac{2a}{\gamma}\Big)}
\end{align*}
In case we want $C_S=\frac{q-1}{2}$, the right hand side of the above inequality should be set to $\frac{1}{q-1}$, i.e, we want $\gamma=2a$. Therefore
\begin{align*}
    B_1=\frac{3\gamma-4a}{\beta}+\Big(2-\frac{2a}{\gamma}\Big)\Big(q-\frac{\gamma}{\beta}\Big)=\frac{2a}{\beta}+\Big(q-\frac{2a}{\beta}\Big)=q>0
\end{align*}
i.e, the condition $B_1 \leq 0$ fails. More generally, we leave the reader to check that if $\frac{a}{\gamma}>0$, all coeffients cannot be nonpositive simultaneously. Therefore $\frac{a}{\gamma} \leq 0$ and the best we can get is that $\lambda \leq \frac{1}{2(q-1)}$. But this is something that is already true for Riemannian manifolds $(M,g)$ with $\text{Ric}(g) \geq g$. This justifies remark 2 and the need for help from the other part of the Hessian, which comes from the next lemma.\newline\newline
\textbf{Lemma 2.4}: For a function $v$ satisfying (*) and $\gamma>0$, $b \in \mathbb{R}$ we have
\begin{align*}
A_2 \int_{M}v^{\gamma-2}|\partial v|^4+B_2\int_{M}v^{\beta+\gamma-\beta q}|\partial v|^2+C_2\int_{M}v^{\gamma}|\partial v|^2+D_2\int_{M}|v_{i\bar{j}}|^2v^{\gamma} \geq 0
\end{align*}
where
\begin{align*}
    & A_2=b^2\Big(1-\frac{1}{n}\Big)+(\beta+1)^2\Big(\frac{2b}{\gamma}-\frac{1}{n}\Big)+2(\beta+1)b\Big(1-\frac{1}{n}\Big),\\
    & B_2=\Big(\frac{2b}{\gamma}-\frac{1}{n}\Big)\Big(q-\frac{\gamma}{\beta}\Big)+\frac{2b}{\beta}\Big(1-\frac{1}{n}\Big),\\
    & C_2=\lambda\Big(\Big(\frac{2b}{\gamma}-\frac{1}{n}\Big)\Big(\frac{\gamma-\beta}{\beta}\Big)-\frac{2b}{\beta}\Big(1-\frac{1}{n}\Big)\Big),\\
    & D_2=\Big(1-\frac{2b}{\gamma}\Big).
\end{align*}
\begin{proof}
    By Cauchy-Schwarz inequality we have
    \begin{align*}
        \int_{M}\Big|v_{i\bar{j}}+b\frac{v_{i}v_{\bar{j}}}{v}\Big|^2 v^{\gamma} \geq \frac{1}{n}\int_{M}\Big(\Box v+b\frac{|\partial v|^2}{v}\Big)^2v^{\gamma}.
    \end{align*}
    Expanding and rearranging terms makes it
    \begin{align*}
&\,\,\,\,\,\,\,\,\,\int_{M}|v_{i\bar{j}}|^2v^{\gamma}+b^2\Big(1-\frac{1}{n}\Big)\int_{M}v^{\gamma-2}|\partial v|^4 +2b\int_{M}v^{\gamma-1}v_{i\bar{j}}v_{\bar{i}}v_{j}\\
        &\,\,\,\,\,\,\,\,-\frac{1}{n}\int_{M}(\Box v)^2v^{\gamma}-\frac{2b}{n}\int_{M}v^{\gamma-1}|\partial v|^2\Box v \geq 0\\
        & \Rightarrow \int_{M}|v_{i\bar{j}}|^2v^{\gamma}+b^2\Big(1-\frac{1}{n}\Big)\int_{M}v^{\gamma-2}|\partial v|^4\\
        &\,\,\,\,\,\,\,\,+2b\Big(\frac{1}{\gamma}\int_{M}(\Box v)^2v^{\gamma}+\int_{M}v^{\gamma-1}|\partial v|^2\Box v-\frac{1}{\gamma}\int_{M}|v_{i\bar{j}}|^2v^{\gamma}\Big)\\
        &\,\,\,\,\,\,\,\,-\frac{1}{n}\int_{M}(\Box v)^2v^{\gamma}-\frac{2b}{n}\int_{M}v^{\gamma-1}|\partial v|^2\Box v \geq 0 \,\,\,\,\text{(By part (3) of lemma 2.2)}\\
        &\Rightarrow \Big(1-\frac{2b}{\gamma}\Big)\int_{M}|v_{i\bar{j}}|^2v^{\gamma}+b^2\Big(1-\frac{1}{n}\Big)\int_{M}v^{\gamma-2}|\partial v|^4+\Big(\frac{2b}{\gamma}-\frac{1}{n}\Big)\int_{M}(\Box v)^2v^{\gamma}\\
        &\,\,\,\,\,\,\,+2b\Big(1-\frac{1}{n}\Big)\int_{M}v^{\gamma-1}|\partial v|^2\Box v \geq 0\\
        & \Rightarrow \Big(1-\frac{2b}{\gamma}\Big)\int_{M}|v_{i\bar{j}}|^2v^{\gamma}+b^2\Big(1-\frac{1}{n}\Big)\int_{M}v^{\gamma-2}|\partial v|^4\\
        &\,\,\,\,\,\,\,+\Big(\frac{2b}{\gamma}-\frac{1}{n}\Big)\Bigg(\frac{\beta q-\gamma}{\beta}\int_{M}v^{\beta+\gamma-\beta q}|\partial v|^2\\
        &\,\,\,\,\,\,\,+\frac{\lambda(\gamma-\beta)}{\beta}\int_{M}v^{\gamma}|\partial v|^2+(\beta+1)^2\int_{M}v^{\gamma-2}|\partial v|^4\Bigg)\\
        &\,\,\,\,\,\,\,+2b\Big(1-\frac{1}{n}\Big)\Big(\frac{1}{\beta}\int_{M}v^{\beta+\gamma-\beta q}|\partial v|^2-\frac{\lambda}{\beta}\int_{M}v^{\gamma}|\partial v|^2+(\beta+1)\int_{M}v^{\gamma-2}|\partial v|^4\Big).\\
    \end{align*}
    Here the last step follows by parts (1) and (2) of lemma 2.2. Collecting like terms gives us the result.
\end{proof}
Now we prove the main theorem:\newline\newline
Consider a parameter $k \geq 0$. By lemmas 2.3 and 2.4 we have
\begin{align*}
    &\int_{M}\Big|v_{\bar{i}\bar{j}}+a\frac{v_{\bar{i}}v_{\bar{j}}}{v}\Big|^2v^{\gamma} \leq (A_1+kA_2)\int_{M}v^{\gamma-2}|\partial v|^4+(B_1+kB_2)\int_{M}v^{\beta+\gamma-\beta q}|\partial v|^2\\
    &\,\,\,\,\,\,\,\,\,\,\,\,\,\,\,\,\,\,\,\,\,\,\,\,\,\,\,\,\,\,\,\,\,\,\,\,\,\,\,\,\,\,\,\,\,\,\,\,\,\,\,\,\,\,\,\,\,\,\,\,\,\,+(C_1+kC_2)\int_{M}v^{\gamma}|\partial v|^2+(D_1+kD_2)\int_{M}|v_{i\bar{j}}|^2v^{\gamma}.
\end{align*}
We note that the right hand side can be rewritten as
\begin{align*}
    & (A_1+kA_2)\int_{M}v^{\gamma-2}|\partial v|^4+(B_1+kB_2)\int_{M}v^{\beta+\gamma-\beta q}|\partial v|^2 +(C_1+kC_2)\int_{M}v^{\gamma}|\partial v|^2\\
    &+(D_1+kD_2)\int_{M}\Big(|v_{i\bar{j}}|^2-\frac{1}{n}(\Box v)^2\Big)v^{\gamma}+\frac{D_1+kD_2}{n}\int_{M}(\Box v)^2v^{\gamma}.
\end{align*}
We replace the last term using part (2) of lemma 2.2. Then we have an inequality (since left hand side of the previous inequality is nonnegative):
\begin{align*}
    A \int_{M} v^{\gamma-2}|\partial v|^4+B \int_{M}v^{\beta+\gamma-\beta q}|\partial v|^2+C\int_{M}v^{\gamma}|\partial v|^2+D \int_{M}\Big(|v_{i\bar{j}}|^2-\frac{1}{n}(\Box v)^2\Big)v^{\gamma} \geq 0
\end{align*}
where
\begin{align*}
    & A=\gamma(\gamma-1)-2a(\gamma-1)+a^2+(3\gamma-4a)(\beta+1)+\Big(2-\frac{2a}{\gamma}\Big)(\beta+1)^2\\
    &\,\,\,\,\,\,\,\,\,\,\,\,\,\,\, + k\Big(b^2\Big(1-\frac{1}{n}\Big)+(\beta+1)^2\Big(\frac{2b}{\gamma}-\frac{1}{n}\Big)+2b(\beta+1)\Big(1-\frac{1}{n}\Big)\Big)+\frac{D}{n}(\beta+1)^2,\\
    & B=\frac{3\gamma-4a}{\beta}+\Big(2-\frac{2a}{\gamma}\Big)\Big(q-\frac{\gamma}{\beta}\Big)+k\Big(\Big(\frac{2b}{\gamma}-\frac{1}{n}\Big)\Big(q-\frac{\gamma}{\beta}\Big)+\frac{2b}{\beta}\Big(1-\frac{1}{n}\Big)\Big)+\frac{D}{n}\Big(q-\frac{\gamma}{\beta}\Big),\\
    & C=\Bigg(\frac{4a-3\gamma}{\beta}+\Big(2-\frac{2a}{\gamma}\Big)\Big(\frac{\gamma}{\beta}-1\Big)+k\Big(\frac{2b}{\gamma}-\frac{1}{n}\Big)\Big(\frac{\gamma}{\beta}-1\Big)\\
    &\,\,\,\,\,\,\,\,\,-\frac{2bk}{\beta}\Big(1-\frac{1}{n}\Big)+\Big(\frac{\gamma}{\beta}-1\Big)\frac{D}{n}\Bigg)\lambda-1,\\
    & D=\frac{2a}{\gamma}-1+k\Big(1-\frac{2b}{\gamma}\Big).
\end{align*}
We want to show that it's possible to choose $\beta \neq 0$, $\gamma$, $a$, $b$ and $k>0$ such that the coefficients $A,B,C,D$ are nonpositive with atleast one of them being negative.\newline
In this direction, we set $D=-\epsilon$ for some $\epsilon \geq 0$, namely by choosing
\begin{align*}
b=\frac{\frac{\gamma\epsilon}{2}+a-\frac{\gamma}{2}}{k}+\frac{\gamma}{2}.
\end{align*}
Then the coefficients become
\begin{align*}
    & A=\gamma(\gamma-1)-2a(\gamma-1)+a^2+(3\gamma-4a)(\beta+1)+(\beta+1)^2\Big(1+\frac{n-1}{n}(k+\epsilon)\Big)\\
    &\,\,\,\,\,\,\,\,\,\,\,\,\,\,\,+\Big(\frac{\frac{\gamma\epsilon}{2}+a-\frac{\gamma}{2}}{k}+\frac{\gamma}{2}\Big)^2k\Big(1-\frac{1}{n}\Big) +2\Big(\frac{\frac{\gamma\epsilon}{2}+a-\frac{\gamma}{2}}{k}+\frac{\gamma}{2}\Big)k(\beta+1)\Big(1-\frac{1}{n}\Big),\\
    &B=\frac{3\gamma-4a}{\beta}+\Big(q-\frac{\gamma}{\beta}\Big)\Big(1+\frac{n-1}{n}(k+\epsilon)\Big)+\frac{2\Big(\frac{\frac{\gamma\epsilon}{2}+a-\frac{\gamma}{2}}{k}+\frac{\gamma}{2}\Big)k}{\beta}\Big(1-\frac{1}{n}\Big),\\
    &C=\Bigg(\frac{4a-3\gamma}{\beta}+\Big(\frac{\gamma}{\beta}-1\Big)\Big(1+\frac{n-1}{n}(k+\epsilon)\Big)-\frac{2\Big(\frac{\frac{\gamma\epsilon}{2}+a-\frac{\gamma}{2}}{k}+\frac{\gamma}{2}\Big)k}{\beta}\Big(1-\frac{1}{n}\Big)\Bigg)\lambda-1,\\
    & D=-\epsilon.
\end{align*}
We observe that these coefficients now are continuous functions of $\gamma$ (in particular continuous at $\gamma=0$). We study the coeffients at $\gamma=0$. In this case $a=bk$ and the coefficients are
\begin{align*}
    & A=2a+a^2+\frac{a^2}{k}\Big(\frac{n-1}{n}\Big)-2a(\beta+1)\Big(\frac{n+1}{n}\Big)+\Big(1+\frac{n-1}{n}(k+\epsilon)\Big)(\beta+1)^2,\\
    &B=-\frac{2a}{\beta}\Big(\frac{n+1}{n}\Big)+q\Big(1+\frac{n-1}{n}(k+\epsilon)\Big),\\
    &C=\Big(\frac{2a}{\beta}\Big(\frac{n+1}{n}\Big)-\Big(1+\frac{n-1}{n}(k+\epsilon)\Big)\Big)\lambda-1,\\
    & D=-\epsilon.
\end{align*}
We set $B=0$, i.e, choose
\begin{align*}
    a=q\Big(1+\frac{n-1}{n}(k+\epsilon)\Big)\frac{n\beta}{2(n+1)}.
\end{align*}
Plugging this value of $a$ in the expression for $A$ and requiring $A \leq 0$ is equivalent to
\begin{align*}
    \beta^2\Big(q^2\Big(1+\frac{n-1}{n}(k+\epsilon)\Big)\frac{(nk+n-1)n}{4k(n+1)^2}-q+1\Big)+\beta\Big(2-\frac{q}{n+1}\Big)+1 \leq 0
\end{align*}
If the discriminant of this quadratic is nonnegative, then existence of a $\beta \neq 0$ satisfying the above inequality is guaranteed. This is equivalent to the condition
\begin{align*}
    &n(n-1)qk^2+k(q(2n^2+n(n-1)\epsilon-2n)-4n^2-4n)\\
    &+q\epsilon(n-1)^2+qn(n-1) \leq 0 \tag{**}
\end{align*}
We look at this expression as a quadratic in $k$. For existence of $k$, it is necessary that the discriminant $\Delta$ is nonnegative. This is equivalent to
\begin{align*}
    \epsilon \leq \frac{4n(n+1)-2q(n-1)-2(n-1)\sqrt{q^2+4qn(n+1)}}{n(n-1)q}
\end{align*}
We now go back to inequality (**). This inequality implies
\begin{align*}
    \frac{4n^2+4n-q(2n^2+n(n-1)\epsilon-2n)+\sqrt{\Delta}}{2n(n-1)q}\geq k \geq \frac{4n^2+4n-q(2n^2+n(n-1)\epsilon-2n)-\sqrt{\Delta}}{2n(n-1)q}.
\end{align*}
Adding $\epsilon$ on both sides to the inequality involving bounds on $k$ we get
\begin{align*}
\frac{4n^2+4n-q(2n^2-n(n-1)\epsilon-2n)+\sqrt{\Delta}}{2n(n-1)q}\geq k+\epsilon \geq \frac{4n^2+4n-q(2n^2-n(n-1)\epsilon-2n)-\sqrt{\Delta}}{2n(n-1)q}.
\end{align*}
We now want to choose $\epsilon>0$ in a way so that the lower bound for $k+\epsilon$ is minimized. This corresponds to a general principle. Consider an expression of the form
\begin{align*}
    \psi(\epsilon)=\alpha\epsilon-\sqrt{\alpha^2\epsilon^2-\beta\epsilon+\gamma}
\end{align*}
with $\alpha,\beta,\gamma>0,\beta^2 > 4\alpha^2\gamma$. Then 
\begin{align*}
    \psi'(\epsilon)=\alpha-\frac{2\alpha^2\epsilon-\beta}{2\sqrt{\alpha^2\epsilon^2-\beta\epsilon+\gamma
    }}>0
\end{align*}
on $\Big(0,\frac{\beta-\sqrt{\beta^2-4\alpha^2\gamma}}{2\alpha^2}\Big)$. Thus $\psi$ is increasing on $\Big(0,\frac{\beta-\sqrt{\beta^2-4\alpha^2\gamma}}{2\alpha^2}\Big)$.
Going back to the lower bound for $k+\epsilon$, this observation means that the lower bound acheives it's minimum for $\epsilon=0$. We have
\begin{align*}
    \Delta_{\epsilon=0}=(4n^2+4n)^2-16(n^2+n)q(n^2-n).
\end{align*}
and the lower bound for $k$ becomes
\begin{align*}
    \frac{2n^2+2n-q(n^2-n)-2\sqrt{(n^2+n)^2-(n^2-n)q(n^2+n)}}{n(n-1)q}.
\end{align*}
With these choices for the parameters $a,b,k,\epsilon$ and $\beta$ we focus our attention on the condition $C<0$. This is equivalent to
\begin{align*}
    \lambda \leq \frac{1}{(q-1)\Big(1+\frac{(n-1)}{n}(k+\epsilon)\Big)}=\frac{1}{(q-1)\Big(\frac{2n+q+2-2\sqrt{(n+1)(n+1-(n-1)q)}}{qn}\Big)}=\frac{1}{2C_S}.
\end{align*}
\section{Proof of theorem B}
Let $u \in C^{\infty}(M)$ be a positive solution of
\begin{align*}
    -\Box u+\lambda u=u^{q}
\end{align*}
We set $u=v^{-\beta}$, $v>0$, $\beta \in \mathbb{R}-\{0\}$ (like we did in the previous section).
Then $v$ satisfies the equation:
\begin{align*}
    \Box v=\frac{1}{\beta}v^{\beta+1-\beta q}-\frac{\lambda}{\beta}v+(\beta+1)\frac{|\partial v|^2}{v} \tag 5
\end{align*}
\textbf{Lemma 3.1}:  For $\gamma \geq 0$, $\gamma \neq \beta q$ we have
\begin{align*}
    &\int_{M}v^{\gamma-1}|\partial v|^2\Box v=\frac{1}{\beta q-\gamma}\int_{M}v^{\gamma}(\Box v)^2-\frac{\lambda(q-1)}{\beta q-\gamma}\int_{M}v^{\gamma}|\partial v|^2\\
    &\,\,\,\,\,\,\,\,\,\,\,\,\,\,\,\,\,\,\,\,\,\,\,\,\,\,\,\,\,\,\,\,\,\,\,\,\,\,\,\,\,\,\,\,\,\,\,\,\,+\frac{(\beta q-\beta-\gamma-1)(\beta+1)}{\beta q-\gamma}\int_{M}v^{\gamma-2}|\partial v|^4.
\end{align*}
\begin{proof}
    We multiply equation (5) by $v^{\gamma}\Box v$ and integrate both sides:
    \begin{align*}
        &\int_{M}v^{\gamma}(\Box v)^2\\
        &=\frac{1}{\beta}\int_{M}v^{\beta+\gamma+1-\beta q}\Box v-\frac{\lambda}{\beta}\int_{M}v^{\gamma+1}\Box v+(\beta+1)\int_{M}v^{\gamma-1}|\partial v|^2\Box v\\
        &=\Big(\frac{\beta q-\beta-\gamma-1}{\beta}\Big)\int_{M}v^{\beta+\gamma-\beta q}|\partial v|^2+\frac{\lambda(\gamma+1)}{\beta}\int_{M}v^{\gamma}|\partial v|^2+(\beta+1)\int_{M}v^{\gamma-1}|\partial v|^2\Box v. \tag 6
    \end{align*}
    where we have applied integration by parts to the first and second terms to obtain the last line.\newline\newline
We multiply (5) by $v^{\gamma-1}|\partial v|^2$:
\begin{align*}
    &\int_{M}v^{\gamma-1}|\partial v|^2\Box v\\
    &=\frac{1}{\beta}\int_{M}v^{\beta+\gamma-\beta q}|\partial v|^2-\frac{\lambda}{\beta}\int_{M}v^{\gamma}|\partial v|^2+(\beta+1)\int_{M}v^{\gamma-2}|\partial v|^4\\
\end{align*}
This implies
\begin{align*}
    \int_{M}v^{\beta+\gamma-\beta q}|\partial v|^2=\beta\int_{M}v^{\gamma-1}|\partial v|^2\Box v+\lambda\int_{M}v^{\gamma}|\partial v|^2-\beta(\beta+1)\int_{M}v^{\gamma-2}|\partial v|^4 \tag 7
\end{align*}
Using equation (7) to eliminate $\int_{M}v^{\beta+\gamma-\beta q}|\partial v|^2$ in equation (6) we get
\begin{align*}
    & \int_{M} v^{\gamma}(\Box v)^2\\
    &=\Big(\frac{\beta q-\beta-\gamma-1}{\beta}\Big)(\beta\int_{M}v^{\gamma-1}|\partial v|^2\Box v+\lambda\int_{M}v^{\gamma}|\partial v|^2-\beta(\beta+1)\int_{M}v^{\gamma-2}|\partial v|^4)\\
    &\,\,\,\,\,\,\,\,\,\,\,+\frac{\lambda(\gamma+1)}{\beta}\int_{M}v^{\gamma}|\partial v|^2+(\beta+1)\int_{M}v^{\gamma-1}|\partial v|^2\Box v\\
    &=(\beta q-\gamma)\int_{M}v^{\gamma-1}|\partial v|^2\Box v+\lambda(q-1)\int_{M}v^{\gamma}|\partial v|^2-(\beta q-\beta-\gamma-1)(\beta+1)\int_{M}v^{\gamma-2}|\partial v|^4.
\end{align*}
This gives us the result.
\end{proof}
\textbf{Lemma 3.2}: For $a\in \mathbb{R}$, $\gamma>0$ we have the inequality:
\begin{align*}
    &(\gamma(\gamma-1)-2a(\gamma-1)+a^2)\int_{M}v^{\gamma-2}|\partial v|^4+(3\gamma-4a)\int_{M}v^{\gamma-1}|\partial v|^2\Box v\\
    &+\Big(2-\frac{2a}{\gamma}\Big)\int_{M}(\Box v)^2v^{\gamma}+\Big(\frac{2a}{\gamma}-1\Big)\int_{M}|v_{i\bar{j}}|^2v^{\gamma}-\int_{M}v^{\gamma}|\partial v|^2 \geq 0
\end{align*}
\begin{proof}
    \begin{align*}
    & 0 \leq \int_{M}\Big|v_{\bar{i}\bar{j}}+a\frac{v_{\bar{i}}v_{\bar{j}}}{v}\Big|^2v^{\gamma}\\
    &=\int_{M}v^{\gamma}|v_{\bar{i}\bar{j}}|^2+a\int_{M}v_{\bar{i}\bar{j}}v_{i}v_{j}v^{\gamma-1}+a\int_{M}v_{ij}v_{\bar{i}}v_{\bar{j}}v^{\gamma-1}+a^2\int_{M}v^{\gamma-2}|\partial v|^4. \tag 8
\end{align*}
We have
\begin{align*}
    & \int_{M}v^{\gamma-1}v_{i}v_{j}v_{\bar{i}\bar{j}}\\
    &=-\int_{M}(v^{\gamma-1}v_{i}v_{j})_{\bar{j}}v_{\bar{i}}\\
    &=-\int_{M}v^{\gamma-1}v_{i\bar{j}}v_{j}v_{\bar{i}}-(\gamma-1)\int_{M}v^{\gamma-2}v_{\bar{j}}v_{i}v_{j}v_{\bar{i}}-\int_{M}v^{\gamma-1}v_{i}\Box v v_{\bar{i}}\\
    &=-\int_{M}v^{\gamma-1}v_{i\bar{j}}v_{\bar{i}}v_{j}-(\gamma-1)\int_{M}v^{\gamma-2}|\partial v|^4-\int_{M}v^{\gamma-1}|\partial v|^2\Box v.
\end{align*}
We note that all the integrals appearing in the last line are real numbers. Therefore the second and third terms in (8) are real and equal. Therefore right hand side of (8) becomes
\begin{align*}
    &\int_{M}v^{\gamma}|v_{\bar{i}\bar{j}}|^2-2a\int_{M}v^{\gamma-1}v_{i\bar{j}}v_{\bar{i}}v_{j}-2a(\gamma-1)\int_{M}v^{\gamma-2}|\partial v|^4\\
    &-2a\int_{M}v^{\gamma-1}|\partial v|^2\Box v +a^2\int_{M}v^{\gamma-2}|\partial v|^4.
\end{align*}
By lemma 2.1 and part (3) of lemma 2.2 this expression is less than or equal to
\begin{align*}
& \gamma(\gamma-1)\int_{M}v^{\gamma-2}|\partial v|^4+\gamma\int_{M}v^{\gamma-1}v_{\bar{i}}v_{j}v_{i\bar{j}}+2\gamma\int_{M}v^{\gamma-1}|\partial v|^2\Box v+\int_{M}v^{\gamma}(\Box v)^2-\int_{M}v^{\gamma}|\partial v|^2\\
&\,\,\,\,\,\,\,\,\,-2a\int_{M}v^{\gamma-1}v_{i\bar{j}}v_{\bar{i}}v_{j}-2a(\gamma-1)\int_{M}v^{\gamma-2}|\partial v|^4-2a\int_{M}v^{\gamma-1}|\partial v|^2\Box v +a^2\int_{M}v^{\gamma-2}|\partial v|^4\\
&=\gamma(\gamma-1)\int_{M}v^{\gamma-2}|\partial v|^4+(\gamma-2a)\Big(\frac{1}{\gamma}\int_{M}(\Box v)^2v^{\gamma}+\int_{M}v^{\gamma-1}|\partial v|^2\Box v-\frac{1}{\gamma}\int_{M}|v_{i\bar{j}}|^2v^{\gamma}\Big)\\
&\,\,\,\,\,\,\,\,\,\,+2(\gamma-a)\int_{M}v^{\gamma-1}|\partial v|^2\Box v +(a^2-2a(\gamma-1))\int_{M}v^{\gamma-2}|\partial v|^4 +\int_{M}v^{\gamma}(\Box v)^2-\int_{M}v^{\gamma}|\partial v|^2\\
    &=(\gamma(\gamma-1)-2a(\gamma-1)+a^2)\int_{M}v^{\gamma-2}|\partial v|^4+(3\gamma-4a)\int_{M}v^{\gamma-1}|\partial v|^2\Box v\\
    &\,\,\,\,\,\,\,\,\,+\Big(2-\frac{2a}{\gamma}\Big)\int_{M}(\Box v)^2v^{\gamma}+\Big(\frac{2a}{\gamma}-1\Big)\int_{M}|v_{i\bar{j}}|^2v^{\gamma}-\int_{M}v^{\gamma}|\partial v|^2.
\end{align*}
\end{proof}
We introduce another parameter in the next lemma.\newline\newline
\textbf{Lemma 3.3}: For $b \in \mathbb{R}$, $\gamma>0$ we have the following inequality:
\begin{align*}
  &\Big(1-\frac{2b}{\gamma}\Big)\int_{M}|v_{i\bar{j}}|^2v^{\gamma}+b^2\Big(1-\frac{1}{n}\Big)\int_{M}v^{\gamma-2}|\partial v|^4+\Big(\frac{2b}{\gamma}-\frac{1}{n}\Big)\int_{M}(\Box v)^2v^{\gamma}\\
        & +2b\Big(1-\frac{1}{n}\Big)\int_{M}v^{\gamma-1}|\partial v|^2\Box v \geq 0\\
\end{align*}
\begin{proof}
        By Cauchy Schwarz inequality we have
    \begin{align*}
        \int_{M}\Big|v_{i\bar{j}}+b\frac{v_{i}v_{\bar{j}}}{v}\Big|^2 v^{\gamma} \geq \frac{1}{n}\int_{M}\Big(\Box v+b\frac{|\partial v|^2}{v}\Big)^2v^{\gamma}.
    \end{align*}
    Expanding and rearranging terms makes it
    \begin{align*}
        &\int_{M}|v_{i\bar{j}}|^2v^{\gamma}+b^2\Big(1-\frac{1}{n}\Big)\int_{M}v^{\gamma-2}|\partial v|^4 +2b\int_{M}v^{\gamma-1}v_{i\bar{j}}v_{\bar{i}}v_{j}\\
        &\,\,\,\,\,\,\,\,\,\,-\frac{1}{n}\int_{M}(\Box v)^2v^{\gamma}-\frac{2b}{n}\int_{M}v^{\gamma-1}|\partial v|^2\Box v \geq 0\\
        & \Rightarrow \int_{M}|v_{i\bar{j}}|^2v^{\gamma}+b^2\Big(1-\frac{1}{n}\Big)\int_{M}v^{\gamma-2}|\partial v|^4 \\
        &\,\,\,\,\,\,\,\,\,+2b\Big(\frac{1}{\gamma}\int_{M}(\Box v)^2v^{\gamma}+\int_{M}v^{\gamma-1}|\partial v|^2\Box v-\frac{1}{\gamma}\int_{M}|v_{i\bar{j}}|^2v^{\gamma}\Big)\\
        &\,\,\,\,\,\,\,\,\,-\frac{1}{n}\int_{M}(\Box v)^2v^{\gamma}-\frac{2b}{n}\int_{M}v^{\gamma-1}|\partial v|^2\Box v \geq 0 \,\,\,\,\text{(By part (3) of lemma 2.2)}\\
        &\Rightarrow \Big(1-\frac{2b}{\gamma}\Big)\int_{M}|v_{i\bar{j}}|^2v^{\gamma}+b^2\Big(1-\frac{1}{n}\Big)\int_{M}v^{\gamma-2}|\partial v|^4+\Big(\frac{2b}{\gamma}-\frac{1}{n}\Big)\int_{M}(\Box v)^2v^{\gamma}\\
        &\,\,\,\,\,\,\,\,\,+2b\Big(1-\frac{1}{n}\Big)\int_{M}v^{\gamma-1}|\partial v|^2\Box v \geq 0\\
\end{align*}
\end{proof}
Combining the inequalities in lemmas (3.2) and (3.3) we get
\begin{align*}
    & \Big(\gamma(\gamma-1)-2a(\gamma-1)+a^2+kb^2\Big(1-\frac{1}{n}\Big)\Big)\int_{M}v^{\gamma-2}|\partial v|^4\\
    &+\Big(2-\frac{2a}{\gamma}+k\Big(\frac{2b}{\gamma}-\frac{1}{n}\Big)\Big)\int_{M}(\Box v)^2v^{\gamma}\\
    &+\Big(3\gamma-4a+2bk\Big(1-\frac{1}{n}\Big)\Big)\int_{M}v^{\gamma-1}|\partial v|^2\Box v\\
    &+\Big(\frac{2a}{\gamma}-1+k\Big(1-\frac{2b}{\gamma}\Big)\Big)\int_{M}|v_{i\bar{j}}|^2v^{\gamma}-\int_{M}v^{\gamma}|\partial v|^2 \geq 0
\end{align*}
for all $k \geq 0$. We now use lemma 3.1 to eliminate the term $\int_{M}v^{\gamma-1}|\partial v|^2\Box v$. The resulting inequality becomes
\begin{align*}
    A(\gamma)\int_{M}v^{\gamma-2}|\partial v|^4+B(\gamma)\int_{M}v^{\gamma}(\Box v)^2+C(\gamma)\int_{M}v^{\gamma}|\partial v|^2+D(\gamma)\int_{M}|v_{i\bar{j}}|^2v^{\gamma} \geq 0
\end{align*}
where
\begin{align*}
    & A(\gamma)=\gamma(\gamma-1)-2a(\gamma-1)+a^2+b^2k\Big(1-\frac{1}{n}\Big)\\
    &\,\,\,\,\,\,\,\,\,\,\,\,\,\,\,\,\,\,\,+\Big(3\gamma-4a+2bk\Big(1-\frac{1}{n}\Big)\Big)\frac{(\beta q-\beta-\gamma-1)(\beta+1)}{\beta q-\gamma},
\end{align*}
\begin{align*}
    B(\gamma)=2-\frac{2a}{\gamma}+k\Big(\frac{2b}{\gamma}-\frac{1}{n}\Big)+\frac{3\gamma-4a+2bk\Big(1-\frac{1}{n}\Big)}{\beta q-\gamma},
\end{align*}
\begin{align*}
    C(\gamma)=-\Bigg(\frac{\lambda(q-1)}{\beta q-\gamma}\Big(3\gamma-4a+2bk\Big(1-\frac{1}{n}\Big)\Big)+1\Bigg),
\end{align*}
\begin{align*}
    D(\gamma)=\Bigg(\frac{2a}{\gamma}-1+k\Big(1-\frac{2b}{\gamma}\Big)\Bigg).
\end{align*}
We choose $b$ such that $D=0$, i.e,
\begin{align*}
    b=\frac{\gamma}{2}\Big(1-\frac{1}{k}\Big)+\frac{a}{k}.
\end{align*}
With this choice of $b$, the coefficients become
\begin{align*}
    &A(\gamma)=\gamma(\gamma-1)-2a(\gamma-1)+a^2+\Big(\frac{\gamma}{2}\Big(1-\frac{1}{k}\Big)+\frac{a}{k}\Big)^2k\Big(1-\frac{1}{n}\Big)\\
    &\,\,\,\,\,\,\,\,\,\,\,\,\,\,\,\,\,\,\,\,\,\,+\Big(3\gamma-4a+2\Big(\frac{\gamma}{2}\Big(1-\frac{1}{k}\Big)+\frac{a}{k}\Big)k\Big(1-\frac{1}{n}\Big)\Big)\frac{(\beta q-\beta-\gamma-1)(\beta+1)}{\beta q-\gamma},\\
    & B(\gamma)=1+\frac{(n-1)k}{n}+\frac{3\gamma-4a+2\Big(\frac{\gamma}{2}\Big(1-\frac{1}{k}\Big)+\frac{a}{k}\Big)k\Big(1-\frac{1}{n}\Big)}{\beta q-\gamma},\\
    &C(\gamma)=-\Bigg(\frac{\lambda(q-1)}{\beta q-\gamma}\Big(3\gamma-4a+2\Big(\frac{\gamma}{2}\Big(1-\frac{1}{k}\Big)+\frac{a}{k}\Big)k\Big(1-\frac{1}{n}\Big)\Big)+1\Bigg).
\end{align*}
We note that now the functions $A,B,C$ are continuous at $\gamma=0$. Plugging in $\gamma=0$ we obtain the inequality:
\begin{align*}
    A\int_{M}v^{-2}|\partial v|^4+B\int_{M}(\Box v)^2+C\int_{M}|\partial v|^2 \geq 0
\end{align*}
where
\begin{align*}
    A=2a+a^2+\frac{a^2}{k}\Big(\frac{n-1}{n}\Big)-2a\Big(\frac{n+1}{n}\Big)\frac{(\beta q-\beta-1)(\beta+1)}{\beta q},
\end{align*}
\begin{align*}
    B=\Big(1+\frac{(n-1)k}{n}-\frac{2a}{\beta q}\Big(\frac{n+1}{n}\Big)\Big),
\end{align*}
\begin{align*}
    C=\Big(\frac{\lambda(q-1)}{\beta q}\Big)\Big(\frac{2a(n+1)}{n}\Big)-1.
\end{align*}
We seek $B \leq 0$ and $A=0$ and set $\frac{a}{\beta}=x$, $\beta=y$. Then
\begin{align*}
    & A=y^2\Big(x+x\Big(\frac{n-1}{kn}\Big)-2\Big(\frac{n+1}{n}\Big)+2\Big(\frac{n+1}{qn}\Big)\Big)\\
    &\,\,\,\,\,\,\,\,\,\,+y\Big(2+2\Big(\frac{n+1}{qn}\Big)-2\Big(\frac{q-1}{q}\Big)\Big(\frac{n+1}{n}\Big)\Big)+2\Big(\frac{n+1}{qn}\Big).
\end{align*}
The existence of $y \neq 0$ such that $A=0$ is equivalent to the discriminant of this quadratic $\Delta$ in $y$ being nonnegative, i.e,
\begin{align*}
    &\Bigg(2+2\Big(\frac{n+1}{qn}\Big)-2\Big(\frac{q-1}{q}\Big)\Big(\frac{n+1}{n}\Big)\Bigg)^2\\
    &-8\Big(\frac{n+1}{qn}\Big)\Bigg(x+x\Big(\frac{n-1}{kn}\Big)-2\Big(\frac{n+1}{n}\Big)+2\Big(\frac{n+1}{qn}\Big)\Bigg) \geq 0
\end{align*}
This is equivalent to the condition
\begin{align*}
    x \leq \frac{(4n^2+4n+q)k}{2(n+1)(kn+n-1)}.
\end{align*}
The condition $B \leq 0$ is equivalent to
\begin{align*}
    x \geq \frac{(kn+n-k)q}{2(n+1)}.
\end{align*}
For these two inequalities to be compatible we must have
\begin{align*}
    & \frac{(kn+n-k)q}{2(n+1)} \leq \frac{(4n^2+4n+q)k}{2(n+1)(kn+n-1)} \\
    & \iff k^2+\Big(2-\frac{4(n+1)}{(n-1)q}\Big)k+1 \leq 0\\
\end{align*}
This is true only if
\begin{align*}
    \frac{2(n+1)}{q(n-1)}-\frac{2(n+1)}{q(n-1)}\sqrt{1-\frac{(n-1)q}{n+1}}-1 \leq k \leq \frac{2(n+1)}{q(n-1)}+\frac{2(n+1)}{q(n-1)}\sqrt{1-\frac{(n-1)q}{n+1}}-1 .
\end{align*}
For $A=0$ and $B<0$ the inequality can be rewritten as
\begin{align*}
    C\int_{M}|\partial v|^2 \geq -B\int_{M}(\Box v)^2 \geq -B\lambda_1\int_{M}|\partial v|^2
\end{align*}
where $\lambda_1$ is the first nonzero eigenvalue of $\Box$ on $(M,\omega)$. If  $v$ is nonconstant this means
\begin{align*}
    &\,\,\,\,\,\,\,\,\,\, C \geq -B\lambda_1\\
    & \Rightarrow \lambda \geq \frac{\lambda_1}{q-1}+\frac{1-\lambda_1\Big(1+\frac{(n-1)k}{n}\Big)}{2\frac{(q-1)(n+1)x}{qn}} \tag{***}
\end{align*}
Our goal now is to try and maximize the right hand side expression of (***). We fix $k$. Then $x$ must lie in the range
\begin{align*}
    \Bigg[\frac{(kn+n-k)q}{2(n+1)},\frac{(4n^2+4n+q)k}{2(n+1)(kn+n-1)}\Bigg].
\end{align*}
We note that the second term on the right hand side of (***) involving $x$ is nonpositive as $\lambda_1\geq 1$. So it achieves it maximum for
\begin{align*}
    x=\frac{(4n^2+4n+q)k}{2(n+1)(kn+n-1)}.
\end{align*}
Plugging this value of $x$, the resulting expression is
\begin{align*}
    F(k)=\frac{1}{q-1}\Bigg(1-\frac{(n+(n-1)k)(kn+n-1)q}{(4n^2+4n+q)k}\Bigg)\lambda_1+\frac{qn(kn+n-1)}{(q-1)(4n^2+4n+q)k}.
\end{align*}
\section{References}
\begin{enumerate}
    \item  M. Bidaut-V\'{e}ron and L. V\'{e}ron. Nonlinear elliptic equations on compact Riemannian manifolds and asymptotics of Emden equations. Invent. Math. 106 (1991), no. 3, 489-539.
    \item   J. Dolbeault, M.J. Esteban, M. Kowalczyk, and M. Loss, Sharp Interpolation Inequalities on the Sphere: New Methods and Consequences, Chinese Annals of Mathematics, Series B, 2013, Volume 34, Issue 1, pp 99-112.
    \item Li, P.; Wang, J. Comparison theorem for K\"{a}hler manifolds and positivity of spectrum, J. Differential Geom. \textbf{69} (2005), no. 1, 43-74.
    \item Chu, J.; Wang, F.; Zhang, K., The rigidity of eigenvalues on K\"{a}hler manifolds with positive Ricci lower bound, J. Reine Angew. Math. \textbf{820} (2025), 213-233.
    \item Datar, V.; Seshadri, H., Diameter rigidity for K\"{a}hler manifolds with positive bisectional curvature, Math. Ann. \textbf{385} (2023), no. 1-2, 471-479.
    \item Baudoin, Fabrice.; Munteanu, Ovidiu., Improved Beckner-Sobolev inequalities on K\"{a}hler manifolds, J. Geom. Anal. \textbf{31} (2021), no. 1, 100-125.
    \item Lott, J. Comparison geometry of holomorphic bisectional curvature for K\"{a}hler manifolds and limit spaces, Duke Math J. \textbf{170} (2021), no. 14, 3039-3071.
    \item Zhang, K. On the optimal volume upper bound for K\"{a}hler manifolds with positive Ricci curvature ( with an appendix by Yuchen Liu), Int. Math. Res. Not. IMRN 2022, no. 8, 6135-6156.
    \item Xiong, Z.; Yang, X. Conjugate radius, volume comparison and rigidity, preprint, arXiv:2408.02080.
    \item Wang, M.; Yang, X. First eigenvalue estimates on complete K\"{a}hler manifolds, preprint, arXiv:2507.09203.
    \item Futaki, A. K\"{a}hler-Einstein metrics and integral invariants, Lecture Notes in Math., 1314 Springer-Verlag, Berlin, 1988. iv+140 pp.
    \item Matsushima Y. Sur la structure du groupe d'hom\'{e}omorphismes analytiques d'unecertaine vari\'{e}t\'{e} k\"{a}hl\'{e}rienne, Nagoya Math. J., \textbf{11} (1957), 145-150.
\end{enumerate}
\end{document}